\documentclass[12pt,a4paper]{article}

\usepackage{amsmath,amssymb}
\usepackage{bm}
\bmdefine{\NNN}{N}

\newcommand{\RRRRR}{{\mathcal R}}
\newcommand{\TTTTT}{{\mathcal T}}

\newcommand{\covered}{\mathrel{<\!\!\!\cdot}}
\newcommand{\ini}{\mathop{\rm in}\nolimits}
\newcommand{\coht}{\mathop{\rm coht}\nolimits}
\newcommand{\height}{\mathrm{ht}}
\newcommand{\tor}{\mathop{\rm Tor}\nolimits}
\newcommand{\define}{\mathrel{:=}}

\newcommand{\gor}{Gorenstein}
\newcommand{\cm}{Cohen-Macaulay}
\newcommand{\type}{{\rm type}}
\newcommand{\qed}{\nolinebreak\rule{.3em}{.6em}}
\newcommand{\joinirred}{join-irreducible}
\newcommand{\rank}{\mathrm{rank}}
\newcommand{\meet}{\wedge}
\newcommand{\join}{\vee}
\newcommand{\lm}{{\rm lm}}
\newcommand{\lt}{{\rm lt}}

\newtheorem{thm}{Theorem}[section]

\newtheorem{lemma}[thm]{Lemma}

\newtheorem{definition}[thm]{Definition}

\newtheorem{claim}{Claim}

\makeatletter
\newcommand{\mysloppy}{\tolerance 9999 \hfuzz .5\p@ \vfuzz .5\p@}
\makeatother
\begin{document}\mysloppy
\begin{center}
\LARGE
A sufficient condition for a Hibi ring to be level 
and levelness of 
Schubert cycles
\end{center}

\begin{center}
\large
Mitsuhiro MIYAZAKI
\\
\normalsize
Dept. Math.,
Kyoto University of Education,
\\
Fukakusa-Fujinomori-cho,
Fushimi-ku, Kyoto, 612-8522 Japan
\\
E-mail:
\tt
g53448@kyokyo-u.ac.jp
\end{center}

\begin{abstract}
Let $K$ be a field, $D$ a finite distributive lattice
and $P$ the set of all join-irreducible elements of $D$.
We show that if 
$\{y\in P\mid y\geq x\}$
is pure for any $x\in P$,
then the Hibi ring $\RRRRR_K(D)$ is level.
Using this result and the argument of sagbi basis theory,
we show that the homogeneous coordinate rings of 
Schubert subvarieties of Grassmannians are level.

\end{abstract}

\section{Introduction}

Let $K$ be a field and $D$ a finite 
distributive lattice.
Hibi \cite{hib} defined the ring $\RRRRR_K(D)$,
which is now called the Hibi ring.
He showed that $\RRRRR_K(D)$ is an 
algebra with straightening law (ASL for short)
over $K$ generated by $D$ and is a normal
affine semigroup ring.
So by the result of Stanley \cite{sta2},
the canonical module of $\RRRRR_K(D)$ is described 
combinatorially.

On the other hand,
Grassmannians and their Schubert subvarieties are
important and fascinating objects of
algebraic geometry and commutative algebra.
And their homogeneous coordinate rings are extensively studied.
For example, it is known that they are normal \cm\ domains
and the homogeneous coordinate rings of Grassmannians are
\gor.
The characterization of the \gor\ property of the 
homogeneous coordinate ring of a Schubert variety
(Schubert cycle for short) is also known.
Note also, that the combinatorial study of 
these rings are the origin of the theory of ASL.

On the other hand, Stanley \cite{sta8} defined the
notion called level for standard graded algebras.
This is a notion between \cm\ property and \gor\ property
for standard graded algebras.

In this article, we give a sufficient condition
for $\RRRRR_K(D)$
to be level.
We also show that the Schubert cycles can be realized as
a subring of a polynomial ring 
which have finite homogeneous sagbi basis.
And the initial algebra of these rings are Hibi rings
satisfying the sufficient condition given in this article.
So the standard deformation argument shows
that every Schubert cycle is level.

\section{Preliminaries}

In this article, all rings and algebras are commutative 
with identity element.

We first recall the results and notation of Hibi \cite{hib},
with a few modication.

We denote by $\NNN$ the set of all non-negative integers.

Let $P$ be a finite partially ordered set (poset for short).

The length of a chain (totally ordered subset) $X$ of $P$
is $\#X-1$, where $\#X$ is the cardinality of $X$.

The rank of $P$, denoted by $\rank P$, is the maximum
of the lengths of chains in $P$.

A poset is said to be
pure if its all maximal chains have the same length.

The height (resp. coheight) of an element $x\in P$,
denoted by $\height_Px$ or simply $\height x$ 
(resp. $\coht_Px$ or $\coht x$),
is the rank of $\{y\in P\mid y\leq x\}$
(resp. $\{y\in P\mid y\geq x\}$).

A poset ideal
 of $P$ is a subset $I$ of
$P$ such that $x\in I$, $y\in P$ and $y\leq x$ 
imply $y\in I$.

For $x$, $y\in P$, $y$ covers $x$, 
denoted by $x\covered y$,
means $x<y$ and 
there is no $z\in P$ such that $x<z<y$.

We denote by $\widehat P$ the extended poset 
$P\cup \{\infty,-\infty\}$ where $\infty$ and
$-\infty$ are new elements and
$-\infty<x<\infty$ for any $x\in P$.

A map $\nu \colon\widehat P\to \NNN$ is order reversing
if $x\leq y$ in $\widehat P$ implies $\nu(x)\geq \nu(y)$
and strictly order reversing if
$x<y$ in $\widehat P$ implies $\nu(x)>\nu(y)$.

The set of all order reversing maps 
(resp. strictly order reversing maps)
from $\widehat P$ to $\NNN$
which map $\infty$ to $0$
is denoted by $\overline\TTTTT(P)$ (resp. $\TTTTT(P)$).

\medskip

Now let $D$ be a finite distributive lattice and $K$ a field.
A \joinirred\ element in $D$ is an element in $D$
which covers exactly one element in $D$.
Recall the result of Birkhoff \cite{bir}.
Let $P$ be the set of all \joinirred\ elements in $D$.
Then $D$ is isomorphic to $J(P)$ ordered by inclusion,
where $J(P)$ is the set of all poset ideals of $P$.
The isomorphisms $\Phi\colon D\to J(P)$ and
$\Psi\colon J(P)\to D$ are given by
\[
\begin{array}{ll}
\Phi(\alpha)\define
	\{x\in P\mid x\leq\alpha\text{ in $D$}\}&\text{for $\alpha\in D$ and}
\\
\Psi(I)\define
	\displaystyle\bigvee_{x\in I}x&\text{for $I\in J(P)$,}
\end{array}
\]
where empty join is defined to be the minimal element of $D$.

Let $\{T_x\}_{x\in \widehat P}$ and $\{X_\alpha\}_{\alpha\in D}$
be families of indeterminates.
Hibi \cite{hib}
defined the ring
$\RRRRR_K(D)$, which is now called the Hibi ring,
by
\[
\textstyle
\RRRRR_K(D)\define
K[\prod_{x\in I\cup\{-\infty\}}T_x\mid I\in J(P)]
\]
as a subring of the polynomial ring
$K[T_x\mid x\in \widehat P]$.
He also showed that $\RRRRR_K(D)$ is isomorphic to
\[
K[X_\alpha\mid \alpha\in D]/
(X_\alpha X_\beta-X_{\alpha\meet\beta}X_{\alpha\join\beta}\mid
\alpha, \beta\in D).
\]
The isomorphism is induced by the $K$-algebra homomorphism
$K[X_\alpha\mid \alpha \in D]\to
K[T_x\mid x\in\widehat P]$
sending $X_\alpha$ to $T_{-\infty}\prod_{x\leq \alpha}T_x$.
And $\RRRRR_K(D)$ is a graded 
ASL
 over $K$
generated by $D$.

It is easily verified that
$\RRRRR_K(D)$ is an affine semigroup ring such that
\[
\RRRRR_K(D)=\bigoplus_{\nu\in\overline\TTTTT(P)}
	K(\prod_{x\in\widehat P}T_x^{\nu(x)}).
\]
So by the result of Hocster \cite{hoc},
$\RRRRR_K(D)$ is a normal \cm\ domain.
And by the result of Stanley \cite{sta2},
\[
\bigoplus_{\nu\in\TTTTT(P)}
	K(\prod_{x\in\widehat P}T_x^{\nu(x)})
\]
is the canonical module of 
$\RRRRR_K(D)$.

\section{A sufficient condtion for $\RRRRR_K(D)$ to be level}

Let $D$ be a finite distributive lattice and 
let $P$ be the set of all \joinirred\ elements in $D$.
We introduce the homogeneous grading on
$\RRRRR_K(D)\simeq
K[X_\alpha\mid \alpha\in D]/
(X_\alpha X_\beta-X_{\alpha\meet\beta}X_{\alpha\join\beta}\mid
\alpha, \beta\in D)$
by setting $\deg X_\alpha=1$ for any $\alpha\in D$.
Then $\RRRRR_K(D)$ is a standard graded algebra,
that is, a Noetherian graded algebra over a field generated
by elements of degree 1.

Stanley \cite{sta8} defined the level property for 
standard graded algebras.
\begin{definition}\rm
Let $A=\bigoplus_{n\geq 0}A_n$ be a \cm{} standard graded 
$K$-algebra,
and let $(h_0,\ldots, h_s)$, $h_s\neq 0$ be the $h$-vector
of $A$.
Then we say that $A$ is level if $h_s=\type A$.
\end{definition}
Suppose that $\dim_K A_1=n$ and $\dim A=d$ in the notation above,
and $S\to A$ is a natural $K$-algebra epimorphism,
where $S$ is a polynomial ring over $K$ with $n$
variables.
Then the following fact is easily verified.

\begin{lemma}\label{level char}
The following conditions are equivalent.
\begin{enumerate}
\item
$A$ is level.
\item
The degree of the generators of the canonical module of
$A$ is constant.
\item
If $F_\bullet$ is the minimal $S$-free resolution
of $A$, then
the degree of the generators of
$F_{n-d}$ is constant.
\end{enumerate}
\end{lemma}

A Hibi ring is not necessarily level as the example
of Hibi \cite[{\S1 e) Example}]{hib} shows.
But we have got the following
sufficient condition for $\RRRRR_K(D)$
to be level.

\begin{thm}\label{suf cond}
Assume that 
$\{y\in P\mid y\geq x\}$ is pure
for any $x\in P$,
then $\RRRRR_K(D)$ is level.
\end{thm}
%
%
Note that the condition on $P$ above is equivalent to the following:
for any $x$, $y\in P$ with $x\covered y$,
$\coht_{\widehat P}x=\coht_{\widehat P}y+1$.

Theorem \ref{suf cond} is a direct consequence of the following 
lemma and the Stanley's description of the canonical
module of a normal affine semigroup ring.

\begin{lemma}\label{suf cond lem}
In the situation of Theorem \ref{suf cond},
for any $\nu\in\TTTTT(P)$
there is $\nu_0\in\TTTTT(P)$ such that
$\nu_0(-\infty)=\rank\widehat P$
and $\nu-\nu_0\in \overline\TTTTT(P)$,
where we set $(\nu-\nu_0)(x)\define\nu(x)-\nu_0(x)$
for any $x\in\widehat P$.
\end{lemma}
{\bf proof}\ \
We set $r=\rank\widehat P$ and
\[
\nu_0(x)\define\max\{\coht_{\widehat P}x,
	r-\nu(-\infty)+\nu(x)\}
\]
for any $x\in \widehat P$.
Then we 
\begin{claim}
For any elements $x$, $y\in\widehat P$ with 
$x\covered y$,
$\nu_0(x)>\nu_0(y)$ and
$\nu(x)-\nu(y)\geq\nu_0(x)-\nu_0(y)$.
\end{claim}
We postpone the proof of the claim and finish
the proof of the lemma
first.

Since 
$\nu_0(x)>\nu_0(y)$ and
$\nu(x)-\nu(y)\geq\nu_0(x)-\nu_0(y)$
for any $x$, $y\in\widehat P$ with $x\covered y$,
it is easily verified that,
$\nu_0(x)>\nu_0(y)$ and
$\nu(x)-\nu(y)\geq\nu_0(x)-\nu_0(y)$
for any $x$, $y\in\widehat P$ with $x<y$.
And $\nu_0(\infty)=0$ by definition.
So $\nu_0\in\TTTTT(P)$.
And
$\nu(x)-\nu_0(x)\geq\nu(y)-\nu_0(y)\geq\nu(\infty)-\nu_0(\infty)=0$
for any $x$, $y\in\widehat P$ with $x<y$.
Therefore
$\nu-\nu_0\in \overline\TTTTT(P)$.
Note $\nu_0(-\infty)=r$ by definition.

Now we prove the claim.
Consider the case where
$\coht_{\widehat P}x\leq
r-\nu(-\infty)+\nu(x)$
first.
Since
$\coht_{\widehat P}y\leq\coht_{\widehat P}x-1$
and $\nu(y)<\nu(x)$, we see that
\[
\coht_{\widehat P}y<\coht_{\widehat P}x\leq r-\nu(-\infty)+\nu(x)
\]
and
\[
r-\nu(-\infty)+\nu(y)<r-\nu(-\infty)+\nu(x).
\]
Therefore
\begin{eqnarray*}
\nu_0(x)&=&r-\nu(-\infty)+\nu(x)\\
	&>&\max\{\coht_{\widehat P}y,r-\nu(-\infty)+\nu(y)\}\\
	&=&\nu_0(y)
\end{eqnarray*}
and
\begin{eqnarray*}
\nu(x)-\nu(y)&=&(r-\nu(-\infty)+\nu(x))-(r-\nu(-\infty)+\nu(y))\\
&\geq&\nu_0(x)-\nu_0(y).
\end{eqnarray*}

Next consider the case where 
$\coht_{\widehat P}x>r-\nu(-\infty)+\nu(x)$.
This happens only when $x\in P$.
Therefore $\coht_{\widehat P}y=\coht_{\widehat P}x-1$
by assumption.
Since
\[
r-\nu(-\infty)+\nu(y)<
r-\nu(-\infty)+\nu(x)\leq\coht_{\widehat P}x-1,
\]
we see that
\[
\coht_{\widehat P}y>r-\nu(-\infty)+\nu(y).
\]
Therefore
\[
\nu_0(y)=
\coht_{\widehat P}y=\coht_{\widehat P}x-1=\nu_0(x)-1<\nu_0(x)
\]
and
\[
\nu_0(x)-\nu_0(y)=1\leq\nu(x)-\nu(y).
\]
So we see the claim.
\qed

We see that the results of Theorem \ref{suf cond} and
Lemma \ref{suf cond lem} are also valid if
$\{y\in P\mid y\leq x\}$ is pure for any $x\in P$,
by using height instead of coheight.

\section{Schubert cycles are level}

Now we fix a field $K$ and integers $m$ and $n$ with
$1\leq m \leq n$.

For an $m\times n$ matrix $M$
with entries in a $K$-algebra $S$,
we denote by $K[M]$ the $K$-subalgebra
of $S$ generated by all the entries of $M$
and by $G(M)$ the $K$-subalgebra of $S$
generated by all maximal minors of $M$.
We also denote by $\Gamma(M)$ the set of all maximal
minors of $M$.

It is known that the homogeneous coordinate ring
of the Grassmann variety $G_m(V)$ of $m$-dimensional
subspaces of an $n$-dimensional $K$-vector space $V$
is $G(X)$,
where $X$ is an $m\times n$ matrix of indeterminates.
It is known that $G(X)$ is an ASL over $K$ generated 
by $\Gamma(X)$,
where we identify $\Gamma(X)$ with a combinatorial object
$\{[c_1,\ldots, c_m]\mid 1\leq c_1<\cdots<c_m\leq n\}$
and define the order of $\Gamma(X)$ by
\[
[c_1,\ldots,c_m]\leq[d_1,\ldots,d_m]
\stackrel{\rm def}{\Longleftrightarrow}
c_i\leq d_i\ \mbox{for $i=1$, \ldots, $m$.}
\]

Let
$0=V_0\subsetneq V_1\subsetneq\cdots\subsetneq V_n=V$
be a complete flag of subspaces of $V$ and let
$a_1$, \ldots, $a_m$ be integers such that
$1\leq a_1<\cdots<a_m\leq n$.
Then the Schubert subvariety $\Omega(a_1,\ldots, a_m)$
of $G_m(V)$ is defined by
\[
\Omega(a_1, \ldots, a_m)\define
\{W\in G_m(V)\mid
\dim(W\cap V_{a_i})\geq i\ \text{for $i=1$, \ldots, $m$}\}.
\]
If we put $b_i=n+1-a_{m+1-i}$ for $i=1$, \ldots, $m$,
$\gamma=[b_1, \ldots, b_m]$ and
$\Gamma(X;\gamma)=\{\delta\in\Gamma(X)\mid
\delta\geq\gamma\}$,
then the homogeneous coordinate ring of the Schubert variety
$\Omega(a_1,\ldots, a_m)$
(Schubert cycle for short) is
\[
G(X;\gamma)\define
G(X)/(\Gamma(X)\setminus\Gamma(X;\gamma))G(X).
\]
This ring is a graded ASL over $K$ generated by $\Gamma(X;\gamma)$
(\cite{dep2}, \cite{bv}).
It is also known that $G(X;\gamma)\simeq
G(U_\gamma)$, where
$U_\gamma$ is the following $m\times n$ matrix 
with independent indeterminates $U_{ij}$ \cite{bv}.
\bgroup
\[
\left(
\begin{array}{cccccccccccc}
0&\cdots&0&U_{1b_1}&\cdots&U_{1b_2-1}&U_{1b_2}&\cdots&\cdots&
U_{1b_m}&\cdots&U_{1n}\\
0&\cdots&0&0&\cdots&0&U_{2b_2}&\cdots&\cdots&
U_{2b_m}&\cdots&U_{2n}\\
\cdots&\cdots&\cdots&\cdots&\cdots&\cdots&\cdots&\cdots&\cdots&
\cdots&\cdots&\cdots\\
0&\cdots&0&0&\cdots&0&0&\cdots&0&U_{mb_m}&\cdots&U_{mn}
\end{array}
\right)
\]
\egroup
We denote the isomorphism
$G(X;\gamma)\stackrel{\sim}{\to}G(U_\gamma)$ by $\Phi$.

Now we introduce a diagonal term order on the polynomial ring
$K[U_\gamma]$.
That is, a monomial order on $K[U_\gamma]$ such that 
the leading monomial of any non-zero minor of $U_\gamma$ is the
product of the main diagonal of it.
For example, the degree lexicographic order induced by
$U_{1b_1}>U_{1b_1+1}>\cdots>U_{1n}>U_{2b_2}>U_{2b_2+1}>\cdots>U_{mn}$.

Then it is easy to see that for standard monomials (of an ASL)
$\mu$ and $\mu'$ on $\Gamma(X;\gamma)$ with $\mu\neq\mu'$,
$\lm(\Phi(\mu))\neq\lm(\Phi(\mu'))$.
So $\Phi(\Gamma(X;\gamma))$ is a sagbi basis of $G(U_\gamma)$,
since $\{\Phi(\mu)\mid\mu$ is a standard monomial on $\Gamma(X;\gamma)\}$
is a $K$-vector space basis of $G(U_\gamma)$.
And 
$\{\lm(\Phi(\mu))\mid\mu$ 
is a standard monomial on $\Gamma(X;\gamma)\}$
is a $K$-vector space basis of the initial subalgebra
$\ini G(U_\gamma)$ of $K[U_\gamma]$.
Note also that
$\lm(\Phi(\alpha))
\lm(\Phi(\beta))=
\lm(\Phi(\alpha\meet\beta))
\lm(\Phi(\alpha\join\beta))$
for any $\alpha$, $\beta\in \Gamma(X;\gamma)$.
So $\ini G(U_\gamma)$ is the Hibi ring 
$\RRRRR_K(\Gamma(X;\gamma))$.

Since the poset of all join-irreducible elements of
$\Gamma(X;\gamma)$ is anti-isomorphic to a finite
poset ideal of $\NNN\times \NNN$ with the componentwise order
\cite{miy},
we see 
by Theorem \ref{suf cond}
that $\ini G(U_\gamma)$ is a level ring.

Now by a standard deformation argument,
we see the following

\begin{lemma}
Let $A$ be a $K$-subalgebra of a polynomial ring
$K[Y_1$, \ldots, $Y_s]$
with
finite homogeneous sagbi basis $f_1$, \ldots, $f_r$.
Let 
$S=K[X_1$, \ldots, $X_r]$
be the polynomial ring with $r$ variables over $K$.
We make $A$ and $\ini A$ $S$-algebras by 
$K$-algebra homomorphisms
$X_i\mapsto f_i$ and
$X_i\mapsto \lt(f_i)$ respectively.
Then
\[
\dim_K\tor_i^S(A,K)_j\leq\dim_K\tor_i^S(\ini A,K)_j
\]
for any $i$ and $j$.
\end{lemma}
Since $A$ and $\ini A$ 
have the same Hilbert function, 
$A$ is \cm, level or \gor\ if so is $\ini A$.
Therefore we see the following
\begin{thm}
Schubert cycles are level.
\end{thm}


\begin{thebibliography}{DEP}
%
%
\bibitem[Bir]{bir}
Birkhoff, G.:
``Lattice Theory.'' Third ed.,
Amer. Math. Soc. Colloq. Publ. No. 25,
Amer. Math. Soc. Providence, R. I. (1967)
%
%
\bibitem[BV]{bv}
Bruns, W. and Vetter, U.:
``Determinantal Rings.''
Lecture Notes in Mathematics {\bf 1327} Springer (1988)
%
%
\bibitem[DEP]{dep2}
DeConcini, C., Eisenbud, D. and Procesi, C.:
``Hodge Algebras.''
Ast\'{e}risque
{\bf 91} (1982)
%
%
\bibitem[Hib]{hib}
Hibi, T.:
{\it Distributive lattices, affine smigroup rings and algebras 
with straightening laws.}
in ``Commutative Algebra and Combinatorics'' (M. Nagata and H. Matsumura, ed.),
Advanced Studies in Pure Math. {\bf11} North-Holland, Amsterdam (1987),
93--109.
%
%
\bibitem[Hoc]{hoc}
Hochster, M.:
{\it Rings of invariants of tori, Cohen-Macaulay rings generated
by monomials and polytopes.}
Ann. of Math. {\bf 96} (1972), 318-337
%
%
%
%
\bibitem[Miy]{miy}
Miyazaki, M.:
{\it On the generating poset of Schubert cycles and
the characterization of
\gor\ property.}
to appear in Bulletin of Kyoto University of Education.
%
%
%
%
\bibitem[Sta1]{sta8}
Stanley R. P.: {\it Cohen-Macaulay Complexes.} in Higher
Combinatorics (M.Aigner, ed.), Reidel, Dordrecht and Boston,
(1977), 51--62.
%
%
\bibitem[Sta2]{sta2}
Stanley, R. P.:
{\it Hilbert Functions of Graded Algebras.}
Adv. Math. {\bf 28} (1978), 57--83.
%
%
%
%
%
%
\end{thebibliography}
\end{document}